\documentclass{amsart}

\numberwithin{equation}{section}

\newtheorem{thm}{Theorem}[section]

\newtheorem{example}[thm]{Example}
\newtheorem{prop}[thm]{Proposition}

\newtheorem{lem}[thm]{Lemma}
\newtheorem{problem}[thm]{Problem}

\usepackage{amssymb}
\begin{document}

\title{On Hahn-Banach type theorems for Hilbert C*-modules}
\author{Michael Frank}
\address{Universit\"at Leipzig, Mathematisches Institut, Augustusplatz 10,
D-04109 Leipzig, Fed.~Rep.~Germany, frank@mathematik.uni-leipzig.de}

\begin{abstract}
We show three Hahn-Banach type extension criteria for (sets of) bounded
C*-linear maps of Hilbert C*-modules to the underlying C*-algebras of
coefficients. One criterion establishes an alternative description of the
property of C*-algebras to be monotone complete or additively complete.
\end{abstract}

\maketitle

\section{Introduction}

The Hahn-Banach theorem belongs to the central results of Banach space theory,
together with the open mapping theorem, the closed graph theorem, and others.
When investigating the generalized Hahn-Banach extension problem for Banach
C*-modules, their C*-submodules and bounded C*-linear maps thereof to the
C*-algebra of coefficients, we must resort to case studies since a general
ans\-wer is out of reach at present. The goal of the present paper is to find
criteria for the particular class of Hilbert C*-modules, i.e.~of Banach
modules over C*-algebras $A$ that admit an $A$-valued inner product,
\cite{Pa1,Lance:95}. We aim to preserve the C*-algebra of coefficients as
the codomain of the extended maps.

\medskip
In comparison to the particular case of Hilbert spaces the problem is
non-trivial since Hilbert C*-submodules of Hilbert C*-modules are in general
not direct summands, and bounded module maps of C*-submodules to the C*-algebra
of coefficients may be non-extendable to bounded module maps on the entire
Hilbert C*-module. To give an example, let $A={\rm C}([0,1])$ denote the
C*-algebra of all continuous functions on the unit interval [0,1] and
$I={\rm C}_0((0,1])$ the subset of all continuous functions vanishing at zero.
The corresponding $A$-valued inner product is defined by $\langle f,g
\rangle_A= fg^*$ for $f,g \in A$. The map $f \in I \longrightarrow g(t)=\sin
(1/t) \cdot f(t) \in I$, \linebreak[4]
$t \in [0,1]$, is a bounded module map on $I$ which cannot
be continued to a bounded $A$-linear map on $A$ preserving the codomain
C*-algebra $A$.
Moreover, the property of a Hilbert C*-submodule $\mathcal N$ of a Hilbert
C*-module $\mathcal M$ of being an orthogonal direct summand thereof sometimes
depends on the choice of the C*-valued inner product on $\mathcal M$, even if
the derived norms are equi\-valent, \cite{Frank:93}. In other words, we are
faced by a difference between the properties of C*-submodules of being an
orthogonal or merely a topological direct summand, cf.~Example \ref{ex44}.

\smallskip
Hahn-Banach extension problems for bounded module maps have several origins
and applications. For Hilbert C*-modules the problem to extend bounded module
maps was first touched in W.~L.~Paschke's papers \cite{Pa1,Pa2}. Knowing
that the analogue of Riesz representation theorem for bounded module functionals
may fail, in general, W.~L.~Paschke analyzed the conditions under which the
$A$-valued inner product on a Hilbert $A$-module $\{ \mathcal M , \langle .,.
\rangle \}$ over a fixed C*-algebra $A$ can be continued to the Banach
$A$-module ${\mathcal M}'$ of all bounded $A$-linear maps $r: \mathcal M \to
A$ turning it into a Hilbert $A$-module, among other similar extension problems.
A particular reinterpretation of this problem is to ask for the existence of
isometric extension of bounded module maps $r$ on $\mathcal M$ to those on
${\mathcal M}'$ preserving the values of $r$ if it is restricted to the
canonically isometrically embedded copy of $\mathcal M$ in ${\mathcal M}'$.
Later on H.~Lin \cite{Lin}, M.~Hamana \cite{Ham92} and the author
\cite{Fr3,FrListe} continued the study of extension problems of bounded module
maps and C*-valued inner products on Hilbert C*-modules and realized the key
role played by the order structure of the positive cone of the underlying
C*-algebra of coefficients. In particular, these results mark the departure
from a topological point of view on Hahn-Banach extension problems for bounded
module maps on pairs of Hilbert C*-modules since there exists a commutative
AW*-algebra (\cite{Kp1}) described by E.~E.~Floyd \cite{Floyd} for which there
does not exist any Frech\'et topology with respect to which both the partial
order structure of its positive cone (which is monotone complete) and its linear
structure are continuous at once. The reader is encouraged to consult the
publications \cite{Ham92,Fr3,FrListe} and, especially, to study H.~Lin's work
\cite{Lin} for more detailed information on the achieved module map extension
results and their proofs. We refer to them at the appropriate places of the
present paper.

A second root and motivation have been results by G.~Vincent-Smith \cite{VS},
G.~Wittstock \cite{Wt2,Wt1}, Ching-Yun Suen \cite{Suen}, P.~S.~Muhly and
Qiyuan Na \cite{MuNa},
A.~M.~Sinclair and R.~R.~Smith \cite{SiSm} on Hahn-Banach type theorems for
completely bounded mo\-dule maps on matricial normed Banach C*-(bi-)modules
(i.e.~operator (bi-)modules over C*-algebras) if the range of the extensions
is allowed to grow inside a (larger) injective C*-algebra, for example the
injective envelope of the C*-algebra of coefficients (\cite{Ham}).
Hilbert C*-modules have the advantage that they can always be represented as
operator modules on Hilbert spaces, and that bounded C*-linear maps between
them are automatically completely bounded, \cite[Prop.~2.7, 2.8]{Wt1}. We
give more details on this point of view in the last section.

\smallskip
To formulate the Hahn-Banach type criteria we obtained the reader should be
aware that a Hilbert $A$-module $\{ \mathcal M, \langle .,. \rangle \}$,
its $A$-dual Banach $A$-module ${\mathcal M}'$ and its $A$-bidual Banach
$A$-module ${\mathcal M}''$ are connected by two canonical isometric embeddings
${\mathcal M} \hookrightarrow {\mathcal M}'' \hookrightarrow {\mathcal M}'$,
where surjectivity may fail at any position and any possible interrelation may
be realized by examples, cf.~\cite{Pa1,Pa2}. Moreover, the process taking
$A$-duals stabilizes since the $A$-dual of ${\mathcal M}''$ equals the
$A$-dual of $\mathcal M$.
Also, the $A$-valued inner product on $\mathcal M$ can be always
continued to an $A$-valued inner product on ${\mathcal M}''$,
cf.~\cite[Th.~2.4]{Pa2}. To define some properties of C*-algebras that come into
the play we say that a C*-algebra $A$ is monotone complete if every bounded
increasingly directed net $ \{ a_\alpha :\alpha \in I \} $ of self-adjoint
elements of $A$ has a least upper bound $ a= \sup \{ a_\alpha : \alpha \in I
\}$ inside $A$. It is said to be additively complete if every norm-bounded
infinite sum of posi\-tive elements possesses a least upper bound inside $A$.
For more details we refer to the next section were the preliminaries will be
explained.

\begin{thm}  \label{th}
   For C*-algebras $A$ the following conditions are equivalent:

   \smallskip
    (i) $\,\,$ for any pair $\{ \{ {\mathcal M}, \langle .,. \rangle \},
      {\mathcal N} \subseteq {\mathcal M} \}$ of a Hilbert $A$-module and a
      Hilbert $A$-sub\-mo\-dule there exists an $A$-linear isometric embedding
      $\phi$ of the Banach $A$-module ${\mathcal N}'$ into the Banach
      $A$-module ${\mathcal M}'$ extending the canonical embedding of $\mathcal
      N$ into ${\mathcal M}'$ induced by the $A$-valued inner product.

   \smallskip
   (ii) $\,\,$ for any pair $\{ \{ {\mathcal M}, \langle .,. \rangle \},
     {\mathcal N} \subseteq {\mathcal M} \}$ of a Hilbert $A$-module and a
     Hilbert $A$-submodule every bounded $A$-linear mapping $r: {\mathcal N}
     \rightarrow A$ can be continued to a bounded $A$-linear mapping $r':
     {\mathcal M} \rightarrow A$ so that
     $\:$(a) $\|r\|=\|r'\|$,
     $\:$(b) $r'$ restricted to $\mathcal N$ equals $r$ and
     $\:$(c) $r'$ restricted to ${\mathcal N}^\bot \subseteq {\mathcal M}$ equals
     zero.

   \smallskip
   (iii) $\,\,$ the multiplier C*-algebra $M(A)$ of $A$ is monotone
     complete.

   \smallskip
   (iv) $\,\,$ the multiplier C*-algebra $M(A)$ of $A$ is additively
     complete.

   \smallskip
   In this case the C*-algebra $A$ has the property that the biorthogonal
   complements ${\mathcal N}^{\bot\bot}$ of Hilbert $A$-submodules $\mathcal N$
   of an arbitrary Hilbert $A$-module $\mathcal M$ estimated with respect to the
   $A$-bidual Hilbert $A$-module ${\mathcal M}''$ of $\mathcal M$ are orthogonal
   direct summands of ${\mathcal M}''$ automatically.
\end{thm}

The implications (iii)$\to$(ii) and (iii)$\to$(i) have been established by
Huaxin Lin \cite[Th.~3.8]{Lin}. The equivalence (iii)$\leftrightarrow$(iv)
has been discovered by K.~Sait{\^o} and J.~D.~M.~Wright \cite[\S 3]{SW}.
The goal of the present paper is the equivalence (i)$\leftrightarrow$(ii) and
the implication (i)/(ii)$\to$(iii)/(iv). 
During the course of proof, we obtain a result valid for particular pairs of
some Hilbert C*-module and one of its Hilbert C*-submodules without any
restrictions to the properties of the C*-algebra of coefficients:

\begin{thm}  \label{th2}
  Let $A$ be a C*-algebra, $\{ {\mathcal M},\langle .,. \rangle \}$ be a Hilbert
  $A$-module and                                            \linebreak[4]
  $\mathcal N \subseteq \mathcal M$ one of its Hilbert $A$-submodules.

  Then every bounded $A$-linear mapping $r: {\mathcal N} \rightarrow A$ can be
  continued to a bounded $A$-linear mapping $r': {\mathcal M} \rightarrow A$
  so that
  $\:$(a) $\|r\|= \|r'\|$,
  $\:$(b) $\, r'$ restricted to $\mathcal N$ equals $r$ and
  $\:$(c) the extended mappings of ${\mathcal N}'$ form a Banach $A$-submodule
   of ${\mathcal M}'$, and
  $\:$(d) the extensions $\{ r'_n : n \in {\mathcal N} \}$ of
   the standardly embedded mappings $\{ r_n = \langle .,n \rangle : n \in
   {\mathcal N} \} \subseteq {\mathcal M}'$ coincide with the latter, $\,$
                      if and only if $\,$
  every bounded $A$-linear mapping $r: {\mathcal N} \rightarrow A$ can be
  continued to a bounded $A$-linear mapping $r': {\mathcal M} \rightarrow A$ so
  that
  $\:$(a) $\|r\|=\|r'\|$,
  $\:$(b) $\, r'$ restricted to $\mathcal N$ equals $r$ and
  $\:$(c) $r'$ restricted to ${\mathcal N}^\bot \subseteq {\mathcal M}$ equals
   zero.

  This happens if and only if the $A$-bidual Hilbert $A$-module ${\mathcal M}''$
  of $\mathcal M$ is the orthogonal direct sum of the orthogonal and the
  biorthogonal complement of $\mathcal N$ inside   \linebreak[4]
  ${\mathcal M}''$ with respect
  to the continued $A$-valued inner product $\langle .,. \rangle$,
  i.e.                                       \linebreak[4]
  ${\mathcal M}''= {\mathcal N}^{\bot\bot} \oplus {\mathcal N}^\bot$.
  Moreover, the $A$-bidual Banach $A$-modules ${\mathcal N}''$ and
  $({\mathcal N}^{\bot\bot})''$ coincide isometrically.

\end{thm}

To weaken the additional assumption on the embedding $\phi: {\mathcal N}'
\to {\mathcal M}'$ to be isometric we can suppose that it has to be merely
bounded and $A$-linear with closed range, cf.~Example \ref{ex44}. Going with
these assumptions we can formulate the following result:

\begin{thm} \label{th3}
  Let $A$ be a C*-algebra, $\{ {\mathcal M}, \langle .,. \rangle \}$ be an
  $A$-reflexive Hilbert $A$-module, i.e.~${\mathcal M} = {\mathcal M}''$.
  Then for Hilbert $A$-submodules ${\mathcal N} \subseteq {\mathcal M}$ the
  following two conditions are equivalent:

  \smallskip
  (i) $\:$ ${\mathcal N}$ is a topological direct summand of ${\mathcal M}$
      (i.e.~not necessarily an orthogonal direct summand).

  \smallskip
  (ii) $\,$ ${\mathcal N} \equiv {\mathcal N}^{\bot\bot} \subseteq
     {\mathcal M}$ and there exists a bounded $A$-linear injective mapping
     $\phi : {\mathcal N}' \to {\mathcal M}'$ such that $\phi(r)[n]=r(n)$ for
     every $r \in {\mathcal N}'$, every $n \in {\mathcal N}$.

  \smallskip
  Without the assumption of $\mathcal M$ being $A$-reflexive (ii) does not
  imply (i), in general.
\end{thm}

\smallskip
The next section contains some definitions and facts from the literature that
have to be stated for proving. The third section is concerned with metric
aspects of the proof, whereas the forth section deals with its inner product
aspects. The proof of the theorems is divided into a number of propositions
which are of independent interest. We close the considerations with a remark
on the structure of particular extensions of C*-valued functionals on Hilbert
C*-modules and on the relation of our results to Hahn-Banach type problems
for operator modules.


\section{{Preliminaries}}

In this section we give definitions and basic facts about Hilbert C*-modules
and C*-algebras necessary for the proofs of the theorems. The papers
\cite{Pa1,Kas,Fr1,Lin:90/2,Lin,Frank:93}, some chapters in \cite{JT,NEWO}, and
the books by E.~C.~Lance \cite{Lance:95} and by I.~Raeburn and D.~P.~Williams
\cite{RaeWil} are used as standard sources of
reference on Hilbert C*-module theory. We make the convention that all
C*-modules of the present paper are left modules by definition. A {\it
pre-Hilbert $A$-module over a C*-algebra} $A$ is an $A$-module $\mathcal M$
equipped with an $A$-valued mapping $\langle .,. \rangle : {\mathcal M}
\times {\mathcal M} \rightarrow A$ which is $A$-linear in the first argument
and has the properties:
\[
\langle x,y \rangle = \langle y,x \rangle^* \; \, , \: \;
\langle x,x \rangle \geq 0 \quad {\rm with} \: {\rm equality} \: {\rm iff}
\quad x=0 \, .
\]
The mapping $\langle .,. \rangle$ is said to be {\it the $A$-valued inner
product on} $\mathcal M$. A pre-Hilbert $A$-module $\{ \mathcal M, \langle
.,. \rangle \}$ is {\it Hilbert} if and only if it is complete with respect
to the norm $\| . \| = \| \langle .,. \rangle \|^{1/2}_A$. We always assume
that the linear structures of $A$ and $\mathcal M$ are compatible. Two Hilbert
$A$-modules are {\it isomorphic} if they are isometrically isomorphic as
Banach $A$-modules, if and only if they are unitarily isomorphic,
\cite{Lance:95}. Moreover, Banach $A$-modules can carry unitarily
non-isomorphic $A$-valued inner products that induce equivalent norms to the
given one, nevertheless, \cite{Frank:93}.

A Hilbert $A$-module $\mathcal M$ is {\it full} if the norm-closed linear span
$\langle {\mathcal M},{\mathcal M} \rangle$ of the values of the $A$-valued
inner product inside $A$ coincides with $A$. A Hilbert $A$-module $\{ \mathcal
M , \langle .,. \rangle \}$ over a C*-algebra $A$ is said to be {\it self-dual}
if and only if every bounded module map $ r: \mathcal M \to A$ is of the form
$\langle . , x_r \rangle$ for some element $x_r \in \mathcal M$. The set of
all bounded module maps $r: \mathcal M \longrightarrow A$ forms a Banach
$A$-module ${\mathcal M}'$. A Hilbert $A$-module is called {\it C*-reflexive}
(or more precisely, {\it $A$-reflexive}) if and only if the map $\Omega$
defined by the formula $\Omega (x)[r]=r(x)$ for each $x \in {\mathcal M}$,
every $r \in {\mathcal M}'$, is a surjective module mapping of $\mathcal M$
onto the Banach $A$-module ${\mathcal M}''$, where ${\mathcal M}''$ consists
of all bounded module maps from ${\mathcal M}'$ to $A$. Note that the property
of being self-dual does not depend on the choice of the C*-algebra of
coefficients $A$ within $\langle {\mathcal M},{\mathcal M} \rangle \subseteq
A \subseteq M(\langle {\mathcal M},{\mathcal M} \rangle)$, whereas the
property of being $A$-reflexive sometimes does. As an example consider
the C*-algebra $A=c_0$ of all sequences converging to zero and set $\mathcal M
=c_0$ with the standard $A$-valued inner product. The multiplier C*-algebra
of $A=c_0$ is $M(A)=l_\infty$, the set of all bounded sequences.
Then ${\mathcal M}'$ equals $l_\infty$ as a one-sided $A$-module, whereas
${\mathcal M}'' = c_0$ again. In contrast, the set of all
bounded $M(A)$-linear maps of ${\mathcal M}'$ to $M(A)$ can be identified
with $l_\infty$. Generally speaking, the $A$-dual Banach $A$-module
${\mathcal M}'$ of a Hilbert $A$-module $\mathcal M$ can be described as
the linear hull of the completed with respect to the topology $\{ \|\langle
x,. \rangle \|_A : x \in {\mathcal M}, \, \|x\| \leq 1 \}$ unit ball of
$\mathcal M$, \cite[Th.~6.4]{Frank:93}.

During our considerations we need a decomposition property of elements of
Hilbert $A$-modules $\mathcal M$ over C*-algebras $A$ that is sharper than
any density argument of $A \circ \mathcal M$ in $\mathcal M$ usually given,
cf.~\cite[p.~5]{Lance:95}. It could be attributed to a number of mathematicians
like E.~Hewitt, B.~E.~Johnson, A.~M.~Sinclair, G.~Skandalis, and others.
We cite a version which may be found in \cite[Th.~4.1]{Ped97} for example.

\begin{lem} \label{decompose}
  Let $A$ be a C*-algebra acting on a Banach $A$-module $\mathcal M$ in such
  a way that $A \circ \mathcal M$ is norm-dense in $\mathcal M$. Then every
  element $x \in \mathcal M$ can be decomposed as $x=ay$ for some $y \in
  \mathcal M$ and some positive $a \in A$ with $\|a\| \leq 1$. Moreover, for
  any given $\varepsilon > 0$ an estimate $\| x-y \| < \varepsilon$ can be
  observed.

  If $\mathcal M$ is a Hilbert $A$-module then by construction the element $a$
  can be selected from the smallest norm-closed two-sided ideal of $A$
  containing the value $\langle x,x \rangle$.
\end{lem}

As an immediate consequence we obtain that any bounded $A$-linear or
$M(A)$-linear map from $\mathcal M$ into $M(A)$ has a range inside $A
\subseteq M(A)$. At the same time every $A$-linear map is automatically
$M(A)$-linear since $M(A)$ is the strict completion of $A$ and the considered
maps are continuous. The $A$-valued inner product
on a Hilbert $A$-module $\mathcal M$ can always be lifted to an $A$-valued
inner product on its $A$-bidual Banach $A$-module ${\mathcal M}''$,
\cite[Th.~2.4]{Pa2}. We already mentioned the canonical isometric embeddings
${\mathcal M} \hookrightarrow {\mathcal M}'' \hookrightarrow {\mathcal M}'$
arising from the C*-valued inner product, and the coincidence $({\mathcal M})'
\equiv ({\mathcal M}'')'$, \cite{Pa2}. A Hilbert $A$-module $\mathcal M$ is
said to be {\it orthogonally complementary as a Hilbert $A$-module} if each
Hilbert $A$-module $\mathcal L$ admitting a bicontinuous embedding of
$\mathcal M$ as a Banach $A$-submodule is orthogonally decomposable as
${\mathcal L}= {\mathcal M} \oplus {\mathcal M}^{\bot}$.

\begin{prop} \label{prop1} {\rm (cf.~\cite[Th.~2.9]{Lin},
                                \cite[Th.~5.6, 6.5]{Frank:93})}
\newline
  Let $A$ be a C*-algebra and $\{ {\mathcal M}, \langle .,. \rangle \}$
  be a full Hilbert $A$-module. Then every bounded module operator on
  $\mathcal M$ possesses an adjoint bounded module operator if and only if
  $\mathcal M$ is orthogonally complementary. The restriction on $\mathcal M$ of
  being full cannot be avoided.
  In this situation either $\mathcal M$ is self-dual or, firstly, the
  $A$-valued inner product $\langle .,. \rangle$ can be extended to an
  $M(A)$-valued inner product on ${\mathcal M}'$ and, secondly, every
  bounded $A$-linear operator on ${\mathcal M}'$ preserves the Banach
  $A$-submodule ${\mathcal M} \hookrightarrow \{ \langle .,x \rangle : x \in
  {\mathcal M} \} \subseteq {\mathcal M}'$ invariant.
\end{prop}

The next statement characterizes monotone complete and additively complete
C*-algebras by alternative conditions. For example, W*-algebras and commutative
C*-algebras $A={\rm C}(X)$ with stonean compact spaces $X$ are monotone
complete, cf.~\cite{Kp1,Lin,Ham92,Fr3}.

\begin{prop}  \label{prop2} {\rm (\cite[Prop.~4.7]{Fr3},
  \cite[Th.~2.2]{Ham92}, \cite[Lemma 3.7]{Lin}, \cite[\S 3]{SW})}
\newline
  A C*-algebra $A$ is monotone complete if and only if $A$ is additively
  complete, if and only if for every Hilbert $A$-module $ \{ \mathcal M,
  \langle   .,. \rangle \}$ the $A$-valued inner product $\langle .,. \rangle$
  on $\mathcal M$ can be continued to an $A$-valued inner product $\langle .,.
  \rangle_D$ on the $A$-dual Banach $A$-module $\mathcal M{\rm '}$ turning
  $\{ \mathcal M{\rm '},\langle .,. \rangle_D \}$ into a self-dual Hilbert
  $A$-module. \newline
  Moreover, in this situation the equalities $\langle x,y \rangle_D =
  \langle x,y \rangle$, $\langle x,r \rangle_D = r(x)$ are satisfied for
  every $x,y \in {\mathcal M} \hookrightarrow {\mathcal M}'$, every $r \in
  {\mathcal M}'$.
\end{prop}

As a proving technique we need a basic construction to switch from a
given Hilbert $A$-module $\mathcal M$ to a bigger Hilbert $A^{**}$-module
${\mathcal M}^\#$ while preserving many useful properties and guaranteeing the
existence and uniqueness of extended ope\-rators and $A$-($A^{**}$-)valued
inner products, (cf.~H.~Lin \cite[Def.~1.3]{Lin:90/2}, \cite[\S 4]{Pa1}).
For example, the C*-valued inner product on a Hilbert W*-module $\mathcal M$
over a W*-algebra $B$ can always be continued to an C*-valued inner product
on its C*-dual Banach W*-module ${\mathcal M}'$ turning it into a self-dual
Hilbert W*-module, \cite[Th.~3.2]{Pa1}. Therefore, for Hilbert W*-modules we
have ${\mathcal M}'' \equiv {\mathcal M}'$. Moreover, taking the biorthogonal
complement ${\mathcal N}^{\bot\bot}$ of a C*-submodule $\mathcal N$ of a
self-dual Hilbert W*-module $\mathcal M$ and forming its C*-dual Banach
W*-module ${\mathcal N}'$ gives the same subset of $\mathcal M$ under the
canonical indentifications usually made. Also, the C*-linear hull of the
completed with respect to the topology $\{ |f(\langle .,x \rangle)| : x \in
\mathcal N \, , \, \, f \in B_* \}$ unit ball of $\mathcal N$ may be identified
with both ${\mathcal N}^{\bot\bot}$ and ${\mathcal N}'$ inside $\mathcal M$
in this situation, \cite[Th.~3.2]{Fr1}. We refer to \cite{Pa1,Lin:90/2,BDH,Fr1}
for more details.

Let $\{ {\mathcal M},\langle .,. \rangle \}$ be a left pre-Hilbert
$A$-module over a fixed C*-algebra $A$. The algebraic tensor
product $A^{**}\otimes {\mathcal M}$ becomes a left
$A^{**}$-module defining the action of $A^{**}$ on its
elementary tensors by the formula $ab \otimes h = a (b \otimes h)$ for
$a,b \in A^{**}$, $h \in \mathcal M$. Now , setting
\[
    \left[ \sum_i a_i \otimes h_i , \sum_j b_j \otimes g_j \right] =
    \sum_{i,j} a_i \langle h_i, g_j \rangle b_j^*
\]
on finite sums of elementary tensors, a degenerate $A^{**}$-valued inner
pre-product is defined. Factorizing $A^{**} \otimes {\mathcal M}$ by $N= \{ z
\in A^{**} \otimes {\mathcal M} : [z,z]=0 \}$ we obtain a pre-Hilbert
$A^{**}$-module denoted by  ${\mathcal M}^{\#}$ in the sequel. It contains
$\mathcal M$ as an $A$-submodule.
If $\mathcal M$ is Hilbert, then ${\mathcal M}^{\#}$ is Hilbert, and vice versa.
More difficult is the transfer of self-duality. If $\mathcal M$ is
self-dual, then ${\mathcal M}^{\#}$ is self-dual, too, but to demonstrate the
converse is an open problem. Other standard properties like
e.g.~C*-reflexivity cannot be transferred.
Every bounded $A$-linear operator $T$ on $\mathcal M$ possesses a unique
extension to a bounded $A^{**}$-linear operator on ${\mathcal M}^{\#}$
preserving the operator norm. Note that the relation $S_1 \bot S_2$ for
subsets $S_1,S_2 \subseteq \mathcal M$ gives rise to the relation $S_1^{\#}
\bot S_2^{\#}$ for the corresponding subsets of ${\mathcal M}^{\#}$.

\begin{lem} \label{r-extend}
  Let $A$ be a C*-algebra and $\mathcal M, \mathcal N$ be two Hilbert
  $A$-modules. The kernel ${\rm ker}(T)$ of every bounded $A$-linear operator
  $T: \mathcal M \to \mathcal N$ coincides with its biorthogonal complement
  with respect to $\mathcal M$. In particular, the kernel of every bounded
  $A$-linear map $r: \mathcal M \to A$ has this property.
  If for some element $r \in {\mathcal M}'$ its kernel possesses a trivial
  orthogonal complement in $\mathcal M$, then $r$ is the zero mapping.
\end{lem}

\begin{proof}
Assume ${\rm ker}(T) \not= {\rm ker}^{\bot\bot}(T)$ inside of $\mathcal M$.
Form the Hilbert $A$-module $\mathcal L = \mathcal M \oplus \mathcal N$ and
consider $T$ as an operator on $\mathcal L$ with kernel ${\rm ker}(T) \oplus
\mathcal N$. Consider the unique extension of $T$ to a bounded $A^{**}$-linear
operator on the self-dual Hilbert $A^{**}$-module $({\mathcal L}^\#)'$, also
denoted by $T$. Taking advantage of the special properties of Hilbert
W*-modules and of the extension operation $\#$ as described above we
obtain that both ${\rm ker}(T)^\#$ and $({\rm ker}(T)^{\bot\bot})^\#$ are
contained in the biorthogonal complement of ${\rm ker}(T) \subseteq \mathcal M
\hookrightarrow ({\mathcal L}^\#)'$ estimated with respect to
$({\mathcal L}^\#)'$ . What is more, the latter is part of the kernel of the
extended to $({\mathcal L}^\#)'$ operator $T$ since $T$ becomes adjointable
there and admits a polar decomposition, \cite[Prop.~15.3.7]{NEWO} and
\cite{Pa1}. We arrive at a contradiction.

Since Banach C*-modules may admit different C*-valued inner products that
induce equivalent norms to the given one and, nevertheless, are not unitarily
equivalent (cf.~\cite{Frank:93}) one further comment is necessary:
the arguments do not depend on the choice of the $A$-valued inner products
on $\mathcal M$ and $\mathcal N$ within their respective unitary equivalence
classes since on $({\mathcal L}^\#)'$ any two $A^{**}$-valued inner products
are unitarily equivalent by self-duality and the kernel of $T$ is always an
orthogonal direct summand there, cf.~\cite{Pa1}.
\end{proof}

\section{{The metric aspects of the proof}}

This section is devoted to the various Banach C*-modules that are related to
a pair $\{ \mathcal M, \mathcal N \subseteq \mathcal M \}$ consisting of a
Hilbert C*-module $\{ \mathcal M, \langle .,. \rangle \}$ and one of its
Hilbert C*-submodules $\mathcal N$, and to their interrelation. The main
emphasis is put on pairwise containment relations, isomorphisms and direct
sum decompositions. The orthogonality of C*-submodules is merely used for
their definition and for their identification with other C*-submodules.
We start with two simple observations concerning C*-reflexive Hilbert
C*-modules and their properties.

\begin{lem} \label{prop-reflex}
  Let $A$ be a C*-algebra and $\{ {\mathcal M} \langle .,. \rangle \}$ be an
  $A$-reflexive, non-self-dual Hilbert $A$-module. Then $\mathcal M$ is the
  largest Banach $A$-submodule ${\mathcal L} \subseteq {\mathcal M}'$ which is a
  Hilbert $A$-module itself, contains the canonical copy of $\mathcal M$
  in ${\mathcal M}'$ as a Hilbert $A$-submodule and possesses the same
  $A$-dual Banach $A$-module ${\mathcal M}'$.

  That is, if there exists a Banach $A$-submodule ${\mathcal L} \subseteq
  {\mathcal M}'$ which is a Hilbert $A$-module itself, contains the canonical
  copy of $\mathcal M$ in ${\mathcal M}'$ as a Hilbert $A$-submodule and is
  larger than $\mathcal M$, then there exists a bounded $A$-linear mapping
  $r: {\mathcal M} \rightarrow A$ that cannot be continued to a bounded
  $A$-linear mapping $r':{\mathcal L} \rightarrow A$.
\end{lem}

\begin{proof}
Let $\mathcal L$ be a Banach $A$-submodule of ${\mathcal M}'$ with the
properties requested. Consider the map from ${\mathcal L}'$ into ${\mathcal M}'$
defined by the restriction of every element of ${\mathcal L}'$ to the
canonical copy of $\mathcal M$ in ${\mathcal M}'$ which has been supposed to
be also contained in $\mathcal L$. This map defines an isometric embedding of
${\mathcal L}'$ into ${\mathcal M}'$, \cite[Prop.~2.1]{Pa2}.
If this map is onto then $\mathcal L$
has to coincide with the canonically embedded copy of $\mathcal M$ in
${\mathcal M}'$ since $\mathcal M$ was supposed to be C*-reflexive.
\end{proof}

\begin{lem} \label{prop-submod}
  Let $A$ be a C*-algebra, $\{ {\mathcal M}, <.,.> \}$ be an $A$-reflexive
  Hilbert $A$-module and ${\mathcal N} \subseteq {\mathcal M}$ be one of its
  Hilbert $A$-submodules. Then the biorthogonal complement
  ${\mathcal N}^{\bot\bot}$ of $\mathcal N$ with respect to $\mathcal M$ is
  $A$-reflexive, too, and             \linebreak[4]
  ${\mathcal N}^{\bot\bot} \cong ({\mathcal N}^{\bot\bot})''$. In general,
  ${\mathcal N}''$ may be smaller than $({\mathcal N}^{\bot\bot})''$.
\end{lem}

\begin{proof}
For every element $p \in {\mathcal M}'$ we obtain a unique element $r_p \in
({\mathcal N}^{\bot\bot})'$ by restricting $p$ to ${\mathcal N}^{\bot\bot}
\subseteq {\mathcal M}$. The resulting map from ${\mathcal M}'$ into
$({\mathcal N}^{\bot\bot})'$ is a contractive $A$-linear mapping of Banach
$A$-modules. By Lemma \ref{r-extend} the kernel of it consists of all elements
of ${\mathcal M}'$ vanishing on $\mathcal N$.
Now, define a map                                   
$\pi: ({\mathcal N}^{\bot\bot})'' \to {\mathcal M}'' \equiv
\mathcal M$ by the rule $\pi(r)(p) := r(r_p)$ for $r \in
({\mathcal N}^{\bot\bot})''$, $p \in {\mathcal M}'$. It is obviously a
module map, and $\pi(({\mathcal N}^{\bot\bot})'')$ is contained in the
biorthogonal complement of $\mathcal N$ with respect to $\mathcal M$ by
construction. Since $\mathcal M$ is C*-reflexive and ${\mathcal N}^{\bot\bot}$
is a Hilbert C*-submodule of $\mathcal M$ we have the equality
$\| \pi(r) \|_{{\mathcal N}^{\bot\bot}} = \| \pi(r) \|_{{\mathcal M}''}$.
However, $\pi(r)$ acts on ${\mathcal N}^{\bot\bot} \hookrightarrow
{\mathcal M}'$ as $r$ acts on ${\mathcal N}^{\bot\bot} \hookrightarrow
({\mathcal N}^{\bot\bot})'$, and therefore
$\| \pi(r) \|_{{\mathcal N}^{\bot\bot}} = \| r \|_{({\mathcal N}^{\bot\bot})''}$
for any $r \in ({\mathcal N}^{\bot\bot})''$. So the map $\pi$ is isometric.
Moreover, the map $\pi$ is onto since ${\mathcal N}^{\bot\bot} \subseteq
({\mathcal N}^{\bot\bot})''$ is fixed under $\pi$.

The isometric inclusion ${\mathcal N}'' \hookrightarrow
({\mathcal N}^{\bot\bot})'' \equiv {\mathcal N}^{\bot\bot}$ follows from Lemma
\ref{prop-reflex} since ${\mathcal N}' \supseteq ({\mathcal N}^{\bot\bot})'$
by Lemma \ref{decompose}, ${\mathcal N} \subseteq {\mathcal N}^{\bot\bot}
\hookrightarrow {\mathcal N}'$ isometrically and Lemma \ref{prop-reflex} is
valid.
However, $\mathcal N$ need not to be equal to ${\mathcal N}^{\bot\bot}$ in
$\mathcal M$.
To give an example set $A= {\mathcal M}= {\rm C}([0,1])$ and ${\mathcal N}=
{\rm C}_0((0,1])$ to obtain ${\mathcal N}'' \equiv {\mathcal N} \subsetneq
{\mathcal N}^{\bot\bot} \equiv ({\mathcal N}^{\bot\bot})'' \equiv {\mathcal M}
\equiv {\mathcal M}''$.
\end{proof}


\medskip
The next step is to prepare the proofs of the first part of Theorem \ref{th2}
and to demonstrate Theorem \ref{th3}.
We will do this by characterizing some special geometric properties of
isometric $A$-linear embeddings $\phi: {\mathcal N}' \to {\mathcal M}'$ for
pairs consisting of a Hilbert $A$-module $\{ {\mathcal M},\langle .,. \rangle
\}$ and its Hilbert $A$-submodule $\mathcal N \subseteq \mathcal M$ over
arbitrary C*-algebras $A$ (if there exist any such embeddings $\phi$ at all).

\begin{prop}  \label{prop-gHB}
   Let $A$ be a C*-algebra, $\{ {\mathcal M} \langle .,. \rangle \}$ be a Hilbert
   $A$-module and $\mathcal N \subseteq \mathcal M$ be a Hilbert $A$-submodule.
   Suppose there is an isometric $A$-linear embedding $\phi$ of the Banach
   $A$-module ${\mathcal N}'$ into the Banach $A$-module ${\mathcal M}'$ which
   extends the standard isometric $A$-linear embedding ${\rm i}: \mathcal N \to
   {\mathcal M}'$, $({\rm i}(n)=\langle .,n \rangle$ for $n \in \mathcal N)$.
   Then $\phi$ has the property that for every $r \in {\mathcal N}'$ the element
   $\phi(r) \in {\mathcal M}'$ equals zero on ${\mathcal N}^\bot \subseteq {\mathcal M}$.
\end{prop}

\begin{proof}
By Lemma 3.5 of \cite{Lin} the map $\phi$ is uniquely defined by the condition
that it extends the canonical embedding of $\mathcal N$ into ${\mathcal M}'$.
Therefore, if $A$ is a W*-algebra we can refer to \cite[Prop.~3.6]{Pa1} for
the properties of $\phi$, and the desired property can be obtained from there.

In case $A$ is an arbitrary C*-algebra we turn to the bidual situation and use
the construction described in the preliminaries, i.e.~switch from the two
Hilbert $A$-modules $\mathcal N \subseteq \mathcal M$ to their canonically
derived Hilbert $A^{**}$-modules ${\mathcal N}^\# \subseteq {\mathcal M}^\#$.
The isometric $A$-linear embedding ${\rm i}: \mathcal N \hookrightarrow
{\mathcal M}'$ extends to an isometric $A^{**}$-linear embedding ${\rm i}^\#$
of ${\mathcal N}^\#$ into the $A^{**}$-dual Hilbert $A^{**}$-module
$({\mathcal M}^\#)'$ of ${\mathcal M}^\#$, cf.~Proposition \ref{prop2}.
Moreover, the isometric $A$-linear embedding $\phi: {\mathcal N}'
\hookrightarrow  {\mathcal M}'$ extends to an isometric $A^{**}$-linear
embedding $\phi:({\mathcal N}^\#)' \hookrightarrow ({\mathcal M}^\#)'$ since
there is a canonical chain of isometric $A$-linear embeddings $\phi
({\mathcal N}') \hookrightarrow {\mathcal M}' \hookrightarrow
({\mathcal M}^\#)'$ by \cite{Pa1} and, furthermore, every element of
$({\mathcal N}^\#)'$ has a unique realization inside the biorthogonal
complement $\phi({\mathcal N}')^{\bot\bot}$ of the set $\phi({\mathcal N}')$
calculated with respect to the Hilbert $A^{**}$-module $({\mathcal M}^\#)'$.
Now, part one of the proof applies and the restriction of every element
$r' \in \phi({\mathcal N}')$ to ${\mathcal N}^\bot \subseteq \mathcal M$
turns out to be zero.
\end{proof}

Going on, we make significant use of Lemma \ref{prop-submod} to show Theorem
\ref{th3}. The example below ensures that the suppositions of Theorem
\ref{th3} are in general weaker than those of the Theorems \ref{th} and
\ref{th2}.

\begin{example}  \label{ex44} {\rm
   Consider the C*-algebra $A=l_\infty$, its two-sided ideal $I=c_0$ and the
   Hilbert $A$-module ${\mathcal M} = A \oplus I$ equipped with the $A$-valued
   inner product $\langle (a,i),(b,j) \rangle = ab^* + ij^*$ for $(a,i),(b,j)
   \in \mathcal M$. The Hilbert $A$-submodule               \linebreak[4]
   ${\mathcal N} = \{ (i,i) : i \in I \}$
   is a topological direct summand of $\mathcal M$ since $\mathcal M$ can be
   decomposed as ${\mathcal M} =
   {\mathcal N} \stackrel{.}{+} \{ (a,0) : a \in A \}$, but it is not an
   orthogonal direct summand. Moreover, ${\mathcal N} \equiv
   {\mathcal N}^{\bot\bot}$. The $A$-dual Hilbert $A$-module ${\mathcal N}' =
   \{ (a,a) : a \in A \}$ can be boundedly and $A$-linearly embedded into the
   Hilbert $A$-module ${\mathcal M}'= A \oplus A$ by the rule $\phi((a,a))
   (b,i) := 2ia^*$ for $a,b \in A$, $i \in I$. The restriction of $\phi
   ({\mathcal N}') \in {\mathcal M}'$ to ${\mathcal N}^\bot = \{ (-i,i) : i
   \in I \} \subset {\mathcal M}$ is obviously non-zero, which contrasts with
   the assumptions (c) and (d) of Theorem \ref{th2}.  }
\end{example}

{\sc Proof of Theorem \ref{th3}.}
We start with the proof of the implication (i)$\to$(ii). Suppose $\mathcal M$
can be decomposed into the topological direct sum ${\mathcal M} = {\mathcal N}
\stackrel{.}{+} \tilde{{\mathcal N}}$. Since ${\mathcal N}$ is the kernel of
a bounded $A$-linear idempotent operator it coincides with its biorthogonal
complement ${\mathcal N}^{\bot\bot}$ inside $\mathcal M$, cf.~Lemma
\ref{r-extend}. For every $r \in {\mathcal N}'$, every $m \in {\mathcal M}$
the mapping $\phi : {\mathcal N}' \to {\mathcal M}'$ can be defined by the
formula $\phi(r)[m]=\phi(r)[n_m+\tilde{n}_m] := r(n_m)$, so (ii) is fulfilled.

To show the converse conclusion (ii)$\to$(i) consider an element $x \in
{\mathcal M}'' \stackrel{i}{\hookrightarrow} {\mathcal M}'$. The element
$x \in {\mathcal M}''$ applied to ${\mathcal M}'$ defines a bounded $A$-linear
functional on the subset $\phi({\mathcal N}') \subseteq {\mathcal M}'$ that
belongs to ${\mathcal N}''$. By Lemma \ref{prop-submod} the biorthogonal
complement ${\mathcal N}^{\bot\bot}$ of $\mathcal N$ with respect to
${\mathcal M}''$ can be identified with ${\mathcal N}''$ preserving the
canonical $A$-valued inner products. Therefore, there exists an element
\linebreak[4]
$y_x \in {\mathcal N}'' \equiv {\mathcal N}^{\bot\bot} \subseteq {\mathcal M}''$
so that $\phi(r)[y_x]= \phi(r)[x]$ for any $r \in {\mathcal N}'$. That is,
$x = y_x + (x-y_x)$, and the summand $(x-y_x) \in {\mathcal M}''$ vanishes on
$\phi({\mathcal N}') \subseteq {\mathcal M}'$. The map $x \in {\mathcal M}''
\to y_x \in {\mathcal N}''$ is bounded and $A$-linear by the assumed properties
of $\phi$. Consequently, ${\mathcal M}''$ admits a topologically direct sum
decomposition as ${\mathcal M}'' = {\mathcal N}^{\bot\bot} \stackrel{.}{+}
\tilde{{\mathcal N}}$, where $\tilde{{\mathcal N}}$ is the set
$\{ x \in {\mathcal M}'' : \phi(r)[x] = 0 \:\, {\rm for} \:\,
r \in {\mathcal N}' \}$. This is condition (i).

\smallskip
The additional assumption of $\mathcal M$ being $A$-reflexive can in general
not be \linebreak[4] dropped.
V.~M.~Manuilov gave an example of a Hilbert W*-submodule $\mathcal N$ in
a non-self-dual Hilbert W*-module $\mathcal M$ that coincides with its
biorthogonal complement with respect to $\mathcal M$ and, nevertheless, is
not a topological direct summand of $\mathcal M$. We realize that $\mathcal N$
additionally has the property (ii) required by Theorem \ref{th3}.

To be more precise, let $A=l_\infty$ be the W*-algebra of all bounded
complex-valued sequences, $c_0$ be the ideal of all sequences converging to
zero and select the Hilbert $A$-module ${\mathcal M} =
l_2(c_0) \oplus l_2(l_\infty)$ with its standard $A$-valued inner product.
Fix the injective operator $T \in {\rm End}_A(l_2(c_0))$ defined by
$T(\{ a_k \}_k) = \{ 1/k \cdot a_k \}_k$ and take the Hilbert $A$-submodule
\[
    {\mathcal N} = \{ (T(x),x) : x \in l_2(c_0) \} \subset \mathcal M \, .
\]
Note, that $\mathcal N$ coincides with its biorthogonal complement
${\mathcal N}^{\bot\bot}$ with respect to $\mathcal M$. Since ${\mathcal M}''
= {\mathcal M}' = l_2(l_\infty)' \oplus l_2(l_\infty)'$ we have ${\mathcal N}'
=  \{ (T(x),x) : x \in l_2(l_\infty)' \} \hookrightarrow {\mathcal M}'$.
The Hilbert $A$-submodule $\mathcal N$ can be realized as the kernel of a
bounded non-adjointable module operator on $\mathcal M$.

However, $\mathcal N$ is not a topological direct summand of $\mathcal M$,
whereas ${\mathcal N}'$ is an orthogonal direct summand of ${\mathcal M}'$.
To see this, assume the existence of a topologically direct decomposition
${\mathcal M} = {\mathcal N} \stackrel{.}{+} {\mathcal L}$ for some Hilbert
$A$-submodule $\mathcal L$ of $\mathcal M$. The resulting $A$-linear bounded
idempotent operator $P: {\mathcal M} \to \mathcal N$ can be described by two
operators $S_1 \in {\rm End}_A(l_2(c_0))$ and $S_2 \in {\rm End}_A(l_2
(l_\infty),l_2(c_0))$ via the formula
\[
   P : {\mathcal M} \to {\mathcal N} \quad ,
                   \quad P((x,y)) = (TS_1(x)+TS_2(y), S_1(x)+S_2(y)) \, .
\]
Relying on the property $({\rm id}_{{\mathcal M}} - P)^2=({\rm id}_{{\mathcal M}} - P)$
we consider the action of $({\rm id}_{{\mathcal M}} - P)$ on the special subsets
$\{ (0,y) : y \in l_2(l_\infty) \}$ and $\{ (x,0) : x \in l_2(c_0) \}$ of
${\mathcal M}$ and obtain the operator equalities
\begin{eqnarray*}
-TS_2+TS_1TS_2-TS_2({\rm id}-S_2) & = & -TS_2 \, ,\\
{\rm id}-TS_1-TS_1({\rm id}-TS_1)+TS_2S_1 & = & {\rm id}-TS_1 \, .
\end{eqnarray*}
By the injectivity of $T$ they can be transformed to the operator equalities
\[
(S_1T+S_2-{\rm id})S_2 = 0 \quad , \qquad (S_1T+S_2-{\rm id})S_1 = 0
\]
which are valid on $l_2(l_\infty)$ and on $l_2(c_0)$, respectively.
Consequently, the operator
\linebreak[4]
$(S_1T+S_2-{\rm id})$ equals zero on
$l_2(c_0) = {\rm ran}(S_1+S_2)$ and hence, on its $A$-dual Hilbert $A$-module
$l_2(l_\infty)'$ by \cite[Prop.~3.6]{Pa1}.
However, $S_2\not= {\rm id}+S_1T$ on $l_2(l_\infty) \subset l_2(l_\infty)'$
since equality would cause $S_1$ to be unbounded by the range restriction to
$S_2$ and $T$. This contradicts the assumption of the existence of a
bounded $A$-linear idempotent operator $P: {\mathcal M} \to {\mathcal N}$.
\hfill $\qed$

\section{{The inner product aspects of the proof}}

To continue the proofs of Theorem \ref{th} and \ref{th2} we have to take into
account aspects of Hilbert C*-modules related to orthogonality of elements,
and of C*-submodules.
So we leave the sphere of considerations on the Banach C*-module level and
deal with orthogonal direct sum decompositions and orthogonal complements.

First, we show that the statements (i) and (ii) of Theorem \ref{th} are
fulfilled for C*-algebras $A$ that possess a monotone complete multiplier
C*-algebra $M(A)$. In particular, the result reproduces
\cite[Th.~4.1]{Wt1} for the special set of operator C*-modules that consists
of all Hilbert C*-modules. However, the restriction to the Hilbert C*-modules
enables us to enlarge the class of well-behaved (in the sense of a Hahn-Banach
theorem) coefficient C*-algebras beyond the set of injective C*-algebras
because this restriction rules out most of the existing more complicated
operator C*-modules.

\begin{prop}  \label{prop-cont} {\rm (cf.~\cite[Th.~3.8]{Lin}) }
  Let $A$ be a C*-algebra that possesses a monotone complete multiplier
  C*-algebra $M(A)$. Let $\{ {\mathcal M}, \langle .,. \rangle \}$ be a
  Hilbert $A$-module and $\mathcal N \subseteq \mathcal M$ be a Hilbert
  $A$-submodule. Every bounded $A$-linear mapping $r:{\mathcal N} \rightarrow
  A$ can be continued to an $A$-linear bounded mapping $r':{\mathcal M}
  \rightarrow A$ with the same norm and the property $r'({\mathcal N}^\bot)=
  \{ 0 \}$. In fact, such an extension is unique.
\end{prop}

\begin{proof}
If $A$ does not contain an identity, then any Hilbert $A$-module
$\mathcal M$ can be considered as a Hilbert $M(A)$-module preserving the
C*-dual Hilbert C*-module ${\mathcal M}'$ since the set $\{ ax : a \in A, x
\in \mathcal M \}$ is norm-dense inside $\mathcal M$ and $A$ is an essential
two-sided ideal of $M(A)$. Consequently, Theorem 3.8 of \cite{Lin}
applies showing the existence of the claimed extension of $r$ with the
stressed for properties.


\end{proof}

\begin{prop}  \label{prop-oc3}
  Let $A$ be a C*-algebra which possesses a monotone complete multiplier
  C*-algebra $M(A)$. Let $\{ {\mathcal M}, \langle .,. \rangle \}$ be a
  Hilbert $A$-module and $\mathcal N \subseteq \mathcal M$ be a Hilbert
  $A$-submodule.
  Then the biorthogonal complement ${\mathcal N}^{\bot\bot}$ of $\mathcal N$
  inside the $A$-bidual Hilbert $A$-module ${\mathcal M}''$ of $\mathcal M$
  estimated with respect to the continued $A$-valued inner product is an
  orthogonal direct summand of ${\mathcal M}''$.
\end{prop}

\begin{proof}
By Lemma 3.7 of \cite{Lin} the $A$-valued inner product on $\mathcal M$ lifts
to an $M(A)$-valued inner product on the $M(A)$-dual Banach
$M(A)$-module ${\mathcal M}'$ of $\mathcal M$, which is the $A$-dual
Banach $A$-module of $\mathcal M$ at the same time and becomes a
self-dual Hilbert $M(A)$-module this way. The biorthogonal complement
${\mathcal N}^{\bot\bot}$ of $\mathcal N$ (embedded into ${\mathcal M}'$ in
the canonical way) with respect to ${\mathcal M}'$ is a self-dual Hilbert
$M(A)$-submodule and direct summand of ${\mathcal M}'$,
\cite[Th.~4.1]{Fr3} and \cite[Th.~2.8]{Fr1} or \cite[Prop.~3.10]{Lin}. Because
of self-duality, the Hilbert $M(A)$-modules ${\mathcal N}^{\bot\bot}
\subseteq {\mathcal M}'$ and ${\mathcal M}'$ equal the $M(A)$-bidual
Hilbert $M(A)$-modules of $\mathcal N$ and $\mathcal M$, respectively,
and ${\mathcal N}^{\bot\bot} \subseteq {\mathcal M}'$ is a direct summand of
${\mathcal M}'$ by \cite[Prop.~3.10]{Lin} or \cite[Th.~2.8]{Fr1}.

\noindent
If $A=M(A)$ then the monotone completeness of $A$ implies the coincidence
of ${\mathcal L}'$ and ${\mathcal L}''$ for any Hilbert $A$-module $\mathcal L$
by \cite{Fr3} or \cite{Lin}. Additional difficulties appear if $A$ is non-unital
because the $A$-bidual Hilbert $A$-module can be different from the
$M(A)$-bidual Hilbert $M(A)$-module of a given Hilbert $A$-module,
in general. So we describe the $A$-bidual Hilbert $A$-modules of $\mathcal N$
and $\mathcal M$ simply as the norm-closed linear hulls of the sets
$\{ a \cdot x : a \in A \, , \: x \in {\mathcal N}^{\bot\bot} \subseteq
{\mathcal M}' \}$ and $\{ a \cdot x : a \in A \, , \: x \in {\mathcal M}' \}$,
respectively, cf.~Lemma \ref{decompose}.
The property of Hilbert $M(A)$-submodules of ${\mathcal M}'$
of being a direct summand of ${\mathcal M}'$ as well as orthogonal complements
are respected by this selection process. Consequently, the biorthogonal
complement ${\mathcal N}^{\bot\bot} \subseteq {\mathcal M}''$ of $\mathcal N$
inside ${\mathcal M}''$ is a direct summand of ${\mathcal M}''$.
\end{proof}

To show the converse we should give some facts about $A$-reflexive Hilbert
C*-modules over C*-algebras $A$ fulfilling assertion (ii) of Theorem \ref{th}:

\begin{prop}  \label{prop-oc}
  Let $A$ be a C*-algebra satisfying assertion (ii) of Theorem \ref{th}.
  Let $\{ {\mathcal M}, \langle .,. \rangle_{\mathcal M} \}$ be an $A$-reflexive
  Hilbert $A$-module.
  Then $\mathcal M$ is orthogonally complementary and every bounded $A$-linear
  operator on $\mathcal M$ possesses an adjoint operator.
\end{prop}

\begin{proof}
Let $\{ {\mathcal L}, \langle .,. \rangle_{\mathcal L} \}$ be another
Hilbert $A$-module containing $\mathcal M$ as a Banach $A$-sub\-module,
i.e.~housing a bicontinuously embedded copy of the Banach $A$-module
$\mathcal M$. We show that both $\mathcal M$ and its biorthogonal complement
${\mathcal M}^{\bot\bot}$ with respect to $\mathcal L$ possess the same
$A$-dual Banach $A$-module ${\mathcal M}'$. Indeed, consider an element
$r$ of the isometrically embedded copy of ${\mathcal M}'$ inside ${\mathcal L}'$
that exists by Theorem \ref{th},(ii) and has the properties claimed there.
Restrict $r$ to the submo\-dule ${\mathcal M}^{\bot\bot} \subseteq \mathcal L$.
We obtain an element $r'$ of the isometrically embedded copy of
$({\mathcal M}^{\bot\bot})'$ inside ${\mathcal L}'$. The difference $(r-r')$
vanishes on $\mathcal M$ since $r$ and $r'$ coincide there, and it vanishes
on ${\mathcal M}^\bot \subseteq \mathcal L$ by the supposed properties of
the embedded $A$-dual Banach $A$-modules formulated at
Theorem \ref{th},(ii). By Lemma \ref{r-extend} the element $(r-r')$ has to
vanish on ${\mathcal M}^{\bot\bot} \subseteq \mathcal L$, too. This proves
our claim. Investigating the pair of Hilbert $A$-modules $\mathcal M \subseteq
{\mathcal M}^{\bot\bot}$ we know that $\mathcal M$ is $A$-reflexive and that
the Hilbert $A$-modules share the same $A$-dual Banach $A$-module. These facts
together with Lemma \ref{prop-reflex} force them to coincide.

Consider an arbitrary element $x \in \mathcal L$ and the restriction of the
$A$-valued bounded functional $\langle .,x \rangle_{\mathcal L}$ to
$\mathcal M$. By supposition condition (ii) of Theorem \ref{th} holds.
So there exists another $A$-valued bounded functional $r_x \in {\mathcal L}'$
with the properties that $r_x$ restricted to $\mathcal M \subseteq \mathcal L$
coincides with $\langle .,x \rangle_{\mathcal L}$ and that $r_x$ restricted to
the orthogonal complement ${\mathcal M}^\bot$ of $\mathcal M$ inside
$\mathcal L$ equals the zero mapping. We want to find a Hilbert $A$-module
containing both $\mathcal M$ isometrically as a Hilbert $A$-submodule and a
copy of $r_x$. Let us show that the value $\langle r_x,r_x \rangle$ has a
meaning and belongs in fact to $A$.

For this aim consider the self-dual Hilbert $A^{**}$-module
$({\mathcal L}^\#)'$ which is constructed from $\mathcal L$ in the canonical
way described in the preliminaries. The two bounded $A$-linear functionals
$\langle .,x \rangle_{\mathcal L}$ and $r_x$ can be continued to elements of
$({\mathcal L}^\#)'$, cf.~\cite{Pa1,Lin:90/2}. The lifting of the $A$-valued
inner product on $\mathcal L$ to an $A^{**}$-valued inner product on
$({\mathcal L}^\#)'$ allows to define inner product values for any two
elements of ${\mathcal L}' \subseteq ({\mathcal L}^\#)'$, in particular for
$\langle .,x \rangle$ and $r_x$. The $A^{**}$-valued functional $(\langle .,x
\rangle -r_x)$ vanishes on $\mathcal M$ by construction. The extension of
$(\langle .,x \rangle -r_x)$ from $\mathcal L$ to $({\mathcal L}^\#)'$
preserves this property and extends its validity to the biorthogonal
complement of $\mathcal M$ taken with respect to the self-dual Hilbert
$A^{**}$-module $({\mathcal L}^\#)'$, cf.~Lemma \ref{r-extend}. However, by
construction and by self-duality the biorthogonal complement of $\mathcal M$
taken with respect to the self-dual Hilbert $A^{**}$-module $({\mathcal L}^\#)'$
contains $r_x$, \cite{Pa1,Fr1}. We obtain the equality
\[
   0 = \langle r_x,x-r_x \rangle_{({\mathcal L}^\#)'}=
   r_x(x)^*-\langle r_x,r_x \rangle_{({\mathcal L}^\#)'}
\]
applying the functional to the particular element $r_x$. In other words,
the value $\langle r_x,r_x \rangle_{({\mathcal L}^\#)'}= \langle r_x,x
\rangle_{({\mathcal L}^\#)'}= r_x(x)^*$ is an element of $A$. Consequently,
the $A$-valued inner product on $\mathcal M$ can be continued to an
$A$-valued inner product on the Banach $A$-submodule of ${\mathcal M}'$
generated by $\mathcal M$ and $r_x$, reducing the $A^{**}$-valued inner
product on $({\mathcal L}^\#)'$ first to the embedded subset ${\mathcal M}'$
and than to the Banach $A$-module of interest. That way the range of this
inner product reduces to a subset of $A$, cf.~\cite{Pa1,Lin:90/2}.

Now, by Lemma \ref{prop-reflex} the assumption $r_x \not\in {\mathcal M}
\hookrightarrow {\mathcal M}'$ leads to a contradiction. Indeed, in this case
we have constructed a Hilbert $A$-module containing $\mathcal M$ as a proper
Hilbert $A$-submodule and possessing the same $A$-dual Banach $A$-module
${\mathcal M}'$ by construction, which contradicts the supposed $A$-reflexivity
of $\mathcal M$.

Summing up, for every element $x \in \mathcal L$ there exists an element
$r_x \in \mathcal M \subseteq \mathcal L$ so that the orthogonal decomposition
$x=r_x \oplus (x-r_x)$ holds inside $\mathcal L$. Hence, $\mathcal M$ is a
direct summand of $\mathcal L$, and $\mathcal M$ is
orthogonally complementary as a Hilbert $A$-module by definition.
To see that any bounded module operator on $\mathcal M$ admits an adjoint
consider $\mathcal M$ as a Hilbert $\langle \mathcal M,\mathcal M
\rangle$-module and apply Proposition \ref{prop1}.
\end{proof}

\begin{prop}  \label{prop-mc}
  Let $A$ be a C*-algebra with property (ii) of Theorem \ref{th}.
  Then the multiplier C*-algebra $M(A)$ of $A$ is monotone complete.
\end{prop}

\begin{proof}
By \cite[Th.~2.4]{Pa2} and Proposition \ref{prop2} the $M(A)$-valued inner
product of each Hilbert $M(A)$-module can be continued to an $M(A)$-valued
inner product on its $M(A)$-bidual Banach $M(A)$-module. Since we want to
show that our suppositions allow the extension to an $M(A)$-valued inner
product on its $M(A)$-dual Banach $M(A)$-module we can suppose
$M(A)$-reflexivity without loss of generality.

So let $\{ {\mathcal M}, \langle .,. \rangle \}$ be an $M(A)$-reflexive
Hilbert $M(A)$-module and consider the subset ${\mathcal N}=\{ x \in
{\mathcal M} : \langle x,x \rangle \in A \}$. The set $\mathcal N$ is
an $A$-module. Indeed, the condition $\langle x,x \rangle \in A$ allows a
representation of $x \in \mathcal M$ as $x=ay$ for some $a \in A$, $y \in
\mathcal M$ by Lemma \ref{decompose} and, therefore, for $x_1,x_2 \in
\mathcal N$ with $x_i=a_i y_i$ for $a_i \in A$, $y_i \in \mathcal M$, $i=1,2$
we obtain
\begin{eqnarray*}
   \lefteqn{\langle x_1+x_2,x_1+x_2 \rangle =} \\ & = &
          \langle a_1 y_1 + a_2 y_2, a_1 y_1 + a_2 y_2 \rangle \\
   & = &  a_1 \langle y_1,y_1 \rangle a_1^* +
          a_2 \langle y_2,y_1 \rangle a_1^* +
          a_1 \langle y_1,y_2 \rangle a_2^* +
          a_2 \langle y_2,y_2 \rangle a_2^*  \in A \, .
\end{eqnarray*}
Hence, the $A$-module axioms are fulfilled by the rules for C*-valued inner
products and ideals. Furthermore, $\mathcal N$ is complete with respect to
the norm on $\mathcal M$.

We want to demonstrate that the $A$-dual Banach $A$-module ${\mathcal N}'$ of
$\mathcal N$ coincides with the $M(A)$-dual Banach $M(A)$-module
${\mathcal M}'$ of $\mathcal M$. Note that any $M(A)$-linear map $r : \mathcal
N \to M(A)$ takes its values in $A \subseteq M(A)$ by Lemma \ref{decompose}.
First, suppose the restriction of two bounded module mappings $r_1,r_2:
{\mathcal M} \to M(A)$ to $\mathcal N$ gives the same element of
${\mathcal N}'$. Then the Hilbert $M(A)$-submodule
$({\rm ker}(r_1-r_2))^\bot \subseteq {\mathcal M}$ does not contain any
element from $\mathcal N$, i.e.~for $x \in ({\rm ker}(r_1-r_2))^\bot$, $x \not=
0$, we have $\langle x,x \rangle \in M(A) \setminus A$. We obtain a two-sided
ideal $\langle ({\rm ker}(r_1-r_2))^\bot, ({\rm ker}(r_1-r_2))^\bot \rangle
\subseteq M(A)$ possessing a trivial intersection with $A$ inside
$M(A)$. The properties of a multiplier algebra imply $({\rm ker}(r_1-
r_2))^\bot = \{ 0 \}$ and $r_1=r_2$ on $\mathcal M$. Consequently,
${\mathcal M}' \subseteq {\mathcal N}'$.
To deal with the second duals we have to use footnotes to distinguish between
the ranges of the C*-linear functionals on $\mathcal N$. Since
${\mathcal M}_{M(A)}' \subseteq {\mathcal N}_{M(A)}'$ we get
\[
  {\mathcal M} \equiv
  {\mathcal M}_{M(A)}'' \supseteq {\mathcal N}_{M(A)}''
  \supseteq  {\mathcal N}_A'' \supseteq {\mathcal N}   \, .
\]
Since $\mathcal N$ is the largest submodule of $\mathcal M$ taking inner
product values in $A$ and \cite[Th.~2.4]{Pa2} guarantees the same property
for the Hilbert $A$-submodule ${\mathcal N}_A''$ of $\mathcal M$ both they
have to coincide, and $\mathcal N$ turns out to be $A$-reflexive.

Now, the assumption to $A$ to fulfil condition (ii) of Theorem \ref{th}
comes into the play. By Proposition~\ref{prop-oc} $\, \mathcal N$ is
orthogonally complementary as a Hilbert $A$-module. In particular, setting
$\mathcal M = M(A)$ we obtain that $\mathcal N = A$ is ortho\-gonally
complementary as a Hilbert $A$-module, too, and that the $A$-dual Banach
$A$-module of $A$ that coincides with the set of all left multipliers of $A$
in fact can be identified with $M(A)$ by Proposition \ref{prop1}.
Referring again to Proposition \ref{prop1} we establish that the $A$-valued
inner product on $\mathcal N \oplus A$ can be continued to an
$M(A)$-valued inner product on ${\mathcal N}' \oplus M(A) \equiv
{\mathcal M}' \oplus M(A)$ turning ${\mathcal M}' \oplus M(A)$
into a self-dual Hilbert $M(A)$-module. So ${\mathcal M}'$ is a
self-dual Hilbert $M(A)$-module, too.

Summarizing, the $M(A)$-valued inner product of every Hilbert
$M(A)$-module can be canonically continued to an $M(A)$-valued
inner product on its $M(A)$-bidual Banach $M(A)$-module
by \cite[Th.~2.4]{Pa2} and, what is more, on its $M(A)$-dual Hilbert
$M(A)$-module. Referring to Proposition \ref{prop2} we are done.
\end{proof}

\begin{prop} \label{prop-oc4}
   Let $A$ be a C*-algebra, $\{ {\mathcal M},\langle .,. \rangle \}$ be a
   Hilbert $A$-module and $\mathcal N \subseteq \mathcal M$ be its Hilbert
   $A$-submodule. Suppose the pair $\{ \mathcal M,\mathcal N\}$ fulfills the
   suppositions of Proposition \ref{prop-gHB}. Then the $A$-bidual Hilbert
   $A$-module ${\mathcal M}''$ of $\mathcal M$ decomposes
   into a direct orthogonal sum as ${\mathcal M}''={\mathcal N}^{\bot\bot}
   \oplus {\mathcal N}^\bot$, and the $A$-bidual Banach $A$-modules
   ${\mathcal N}''$ and $({\mathcal N}^{\bot\bot})''$ coincide.
   The converse conclusion is also true.
   Both the Hilbert $A$-submodules $\mathcal N$ and ${\mathcal N}^{\bot\bot}
   \subseteq \mathcal M$ are not in general direct summands of $\mathcal M$
   or ${\mathcal M}''$.
\end{prop}

\begin{proof}
The $A$-valued inner product on $\mathcal M$ can be continued to an $A$-valued
inner product on the $A$-bidual Banach $A$-module ${\mathcal M}''$ of $\mathcal
M$ by \cite[Th.~2.4]{Pa2}. Moreover, the biorthogonal complement
${\mathcal N}^{\bot\bot}$ of $\mathcal N$ with respect to ${\mathcal M}''$ is
isomorphic to the $A$-bidual Banach $A$-module
$({\mathcal N}^{\bot\bot})''$ by Lemma \ref{prop-reflex}.

A discussion similar to that in the proof of Proposition \ref{prop-oc}
shows that every element $\langle .,m \rangle \in {\mathcal M}'$ with $m \in
{\mathcal M}''$ can be decomposed into the direct sum of an element $r_m \in
{\mathcal N}^{\bot\bot}$ and an element $(\langle .,m \rangle -r_m) \in
{\mathcal N}^\bot$. Consequently, ${\mathcal M}'' = {\mathcal N}^{\bot\bot}
\oplus {\mathcal N}^\bot$ implying ${\mathcal M}' =
({\mathcal N}^{\bot\bot})' \oplus ({\mathcal N}^\bot)'$ and
$({\mathcal N}^{\bot\bot})' \equiv {\mathcal N}'$ by the supposed property (ii)
of Theorem \ref{th2}. Refering to Lemma \ref{prop-reflex} we obtain the
identity $({\mathcal N}^{\bot\bot})'' \equiv {\mathcal N}''$

The second assertion on the converse conclusion follows from Lemma
\ref{prop-reflex} if we continue bounded module mappings $r \in
({\mathcal N}^{\bot\bot})'$ to bounded module mappings $r' \in
{\mathcal M}'$ by simply setting $r'=r \circ P$, where $P: {\mathcal M}'' \to
{\mathcal N}^{\bot\bot}$ is the projection existing by assumption.

To give an example of a situation where $\mathcal N$ is not a direct summand
of                                             \linebreak[4]
${\mathcal M} \equiv {\mathcal M}''$, but its biorthogonal complement
${\mathcal N}^{\bot\bot}$ is, consider the W*-algebra   \linebreak[4]
$A=\mathcal M$
of all bounded linear operators on a separable Hilbert space as a Hilbert
$A$-module over itself and set $I = \mathcal N$ to be the Hilbert $A$-submodule
of all compact operators. For this pair the generalized Hahn-Banach theorem
holds.

As an example for the case ${\mathcal M} \not\equiv {\mathcal M}''$ consider
the Hilbert $A$-module ${\mathcal M} = A \oplus I$ equipped with its canonical
$A$-valued inner product for the same W*-algebra $A$. The Hilbert $A$-submodule
${\mathcal N} = \{ (i,i) : i \in I \}$ coincides with its biorthogonal
complement taken with respect to $\mathcal M$, but it is not an orthogonal
direct summand of $\mathcal M$. By contrast its biorthogonal complement
${\mathcal N}^{\bot\bot} = \{ (a,a) : a \in A \}$ in the $A$-bidual  Hilbert
$A$-module ${\mathcal M}'' = A \oplus A$ of $\mathcal M$ certainly is an
orthogonal direct summand.
\end{proof}

The following example shows that the bidual Banach C*-modules of a Hilbert
C*-sub\-module and of its biorthogonal complement do not in general coincide:
if we set $A={\mathcal M}={\rm C}([0,1])$ and ${\mathcal N}={\rm C}_0((0,1])$,
then ${\mathcal M} = {\mathcal N}^{\bot\bot} \subseteq \mathcal M$, but the
Hahn-Banach theorem does not hold for the pair $\{ {\mathcal N},{\mathcal M}
\}$ of Hilbert $A$-modules (see introduction) since both
${\mathcal N}^{\bot\bot} \equiv {\mathcal M}$ and $\mathcal N$ are
$A$-reflexive, but different.

\medskip
{\sc Proof of Theorem \ref{th} and \ref{th2}.}
Proposition \ref{prop-gHB} demonstrates that the first set of conditions
stated at Theorem \ref{th2} implies the second one. To show the converse
implication consider two elements $r_1,r_2 \in {\mathcal N}'$, an element
$a \in A$ and the extended to ${\mathcal M}'$ elements $s:=(r_1+r_2)'-r_1'
-r_2'$, $t:=a(r')-(ar)'$. By assumption $s,t \in {\mathcal M}'$ vanish on
both $\mathcal N$ and on ${\mathcal N}^\bot \subseteq \mathcal M$. Since the
orthogonal complement of $\mathcal N \oplus {\mathcal N}^\bot$ with respect
to $\mathcal M$ is trivial Lemma \ref{r-extend} applies, and $s=t=0$ on
$\mathcal M$. So by assumption the subset of ${\mathcal M}'$ consisting of
all extended elements of ${\mathcal N}'$ has the structure of a Banach
$A$-module inherited from that of the Banach $A$-module ${\mathcal N}'$.
Moreover, for $x \in \mathcal N$ the element $(\langle .,x \rangle - (\langle
.,x \rangle)') \in {\mathcal M}'$ vanishes on both $\mathcal N$ and ${\mathcal
N}^\bot$, too, by assumption. The same argument as above yields $\langle .,x
\rangle \equiv (\langle .,x \rangle)'$ for any $x \in \mathcal N$.

By the way we have obtained the equivalence of the assertions (i) and (ii)
of Theorem \ref{th}. Together with the Propositions \ref{prop-cont},
\ref{prop-mc}, \ref{prop-oc3}, \ref{prop-oc4} and \ref{prop2} the proofs of
Theorem \ref{th} and \ref{th2} are complete.
\hfill $\qed$

\section{{Final remarks}}

Hahn-Banach type problems for Hilbert C*-modules are closely related to those
for general operator modules as treated by G.~Vincent-Smith \cite{VS},
G.~Wittstock \cite{Wt2,Wt1}, Ching-Yun Suen \cite{Suen},
P.~S.~Muhly and Qiyuan Na \cite{MuNa},
A.~M.~Sinclair and R.~R.~Smith \cite{SiSm}, and others. Hilbert C*-modules
serve as a special class of examples of operator modules,
\cite[Prop.~2.7, 2.8]{Wt1}. All these Hahn-Banach type theorems for
completely bounded mo\-dule maps on matricial normed Banach C*-(bi-)modules
(i.e.~operator (bi-)modules over C*-algebras) are formulated under the
assumption that the the range of the extensions is allowed to grow inside a
(larger) injective C*-algebra, for example the injective envelope of the
C*-algebra of coefficients (\cite{Ham}). Naturally, order aspects play a
mayor role in these investigations, cf.~\cite{MuNa}. All these theorems give
only sparse information about situations with sharper codomain restrictions,
so Hilbert C*-modules could serve as a model to develop more general
Hahn-Banach type criteria for operator (bi-)modules.
Our presented case study demonstrates that monotone complete C*-algebras
serve as an upper bound for attempts to extend these results to a larger
class of C*-algebras of coefficients and codomains. However, the sharpest
restriction has been obtained in collaboration with V.~I.~Paulsen during a
research stay of the author in Houston in 1998. The example belongs to another
class of operator bimodules: taking an arbitrary C*-algebra $A$ and its
injective envelope $I(A)$, a monotone complete C*-algebra (\cite{Ham}), so
this pair of operator $A$-bimodules fulfills a Hahn-Banach type extension
theorem for completely bounded $A$-linear maps into $A$ if and only if the
C*-algebra $A$ is injective itself (and coincides with its injective envelope).
So the Wittstock extension theorems for operator (bi-)modules turn out to be
actually already criteria. Furthermore, after a talk of the author given at a
meeting on operator spaces at CIRM in Marseille in July 1998 M.~Junge formulated
the idea that for an arbitrary pair of an operator bimodule and one of its
operator subbimodules over a given C*-algebra $A$ any completely bounded module
map of the submodule to $A$ admits a completely bounded extension to the entire
operator bimodule preserving the completely boundedness norm and taking values
at most in the bidual to $A\,$ W*-algebra $A^{**}$, if and only if $A$ has the
weak expectation property in the sense of E.~C.~Lance \cite{Lance}. The proof
was found during discussions of M.~Junge with the author at that meeting.
These results will be published elsewhere since they lead beyond the topic of
the presented paper, see \cite{FrPa:99}, \cite{Fr:2000}.

\smallskip
As a second addition we want to explain another equivalent condition
to the assertions of Theorem~\ref{th} which is much more technical. We get an
idea of how particular mappings $r \in {\mathcal N}'$ may be isometrically
extended to mappings $r' \in {\mathcal M}'$ failing to be zero on
${\mathcal N}^\bot \subseteq \mathcal M$ (cf.~Example \ref{ex44}):

\begin{prop}  \label{prop-sing}
  Let $A$ be a C*-algebra with the property that for every Hilbert $A$-module
  $\mathcal M$ and every Hilbert $A$-submodule $\mathcal N \subseteq \mathcal M$
  every bounded $A$-linear mapping $r: \mathcal N \to A$ can be continued to a
  bounded $A$-linear mapping $r':\mathcal M \to A$ possessing the same norm so
  that the restriction of $r'$ to $\mathcal N$ equals $r$. If additionally
  the equality $\langle r,r \rangle = \langle r',r' \rangle \in A^{**}$ is
  supposed to hold for the $A^{**}$-valued inner product on the self-dual
  Hilbert $A^{**}$-module $({\mathcal M}^\#)'$, then the restriction of each
  extension $r'$ to ${\mathcal N}^\bot \subseteq \mathcal M$ equals zero and
  hence, $A$ has the property (ii) of Theorem \ref{th}.
\end{prop}

\begin{proof}
The mapping $r': {\mathcal M} \to A$ canonically continues to an element of
$({\mathcal M}^\#)'$ and possesses a canonical orthogonal decomposition
$r'=r_1 \oplus r_2$ with $r_1,r_2 \in ({\mathcal M}^\#)'$, $r_1$ restricted to
$({\mathcal N}^\#)^\bot \subseteq ({\mathcal M}^\#)'$ equals zero and $r_2$
restricted to ${\mathcal N}^\# \subseteq ({\mathcal M}^\#)'$ equals zero, too.
The Banach $A$-module ${\mathcal N}'$ is canonically and isometrically embedded
into $({\mathcal N}^\#)' \equiv ({\mathcal N}^\#)^{\bot\bot} \subseteq
({\mathcal M}^\#)'$. Consequently, $r=r_1$ and the equality
\[
\langle r',r' \rangle = \langle r,r \rangle + \langle r_2,r_2 \rangle \in
A^{**}
\]
holds for the appropriate values of the $A^{**}$-valued inner product on
$({\mathcal M}^\#)'$. Hence, $r_2$ has to be zero to fulfill the additional
assumption.
\end{proof}

The point is that the metric equality
$\|r\|=\|r'\|=\|\langle r,r \rangle + \langle r_2,r_2 \rangle \|_A^{1/2}$
can in general be valid for (existing) extensions $r'$ of $r$ even if $r_2
\not=0$.
We have to conclude that the consideration of any possible extensions
$r' \in {\mathcal M}'$ of a given bounded $A$-linear mapping $r: {\mathcal N}
\to A$ (if there is any at all) has to be made separately in every particular
situation.

For example, consider the C*-algebra $A={\rm C}([0,1])$ of all continuous
functions on the unit interval [0,1] as a Hilbert $A$-module $\mathcal M$ over
itself. Set ${\mathcal N} = \{ f \in A : f \equiv 0 \:\, {\rm on} \:\, [1/2,1]
\}$. The embedding map ${\rm i}: {\mathcal N} \to A$ is an element of
${\mathcal N}'$ which does not belong to the subset $\{ \langle .,n \rangle :
n \in {\mathcal N} \} \subset {\mathcal N}'$. Every function $f_i \in A$
taking the value $1$ on the interval [0,1/2] and satisfying the equality
$\|f_i\|=1$ is an extension of ${\rm i} \in {\mathcal N}'$ inside ${\mathcal M}'
=A$. However, ${\mathcal N}^{\bot\bot}$ is not a direct summand of $A=\mathcal
M$ and there does not exist any extension $f_i$ of ${\rm i}$ in $\mathcal M$
whose restriction to ${\mathcal N}^\bot$ equals zero.

\begin{problem}
Denote by $\mathcal M$ and $\mathcal N \subseteq \mathcal M$ a Hilbert
C*-module and one of its Hilbert C*-submodules.
For which kind of C*-algebras $A$ is the biorthogonal complement
${\mathcal N}^{\bot\bot}$ of $\mathcal N$ inside $\mathcal M$ or inside its
$A$-bidual Hilbert $A$-module ${\mathcal M}''$, respectively, always a
topological or orthogonal direct summand?
\end{problem}

\medskip
\textbf{Acknowledgement:} The author thanks D.~P.~Blecher, K.-D.~K\"ursten,
\linebreak[4]
V.~M.~Ma\-nuilov, V.~I.~Paulsen and E.~V.~Troitsky for valuable discussions and
comments.



\end{document}